\documentclass{article}
\usepackage[T2A]{fontenc}
\usepackage[utf8]{inputenc}
\usepackage[english,russian]{babel}
\usepackage[tbtags]{amsmath}
\usepackage{amsfonts,amssymb,mathrsfs,amscd}
\usepackage[hyper]{msb-a}

\makeatletter
\gdef\No{{\select@language{russian}\textnumero}}
\makeatother

\JournalName{}
\numberwithin{equation}{section}

\def \Z {{\mathbf {Z}}}
\def \JJ {{\mathbf {J}}}
\def \N {{\mathbf {N}}}
\def \un {{\bf 1}}
\def \B {{\cal B}}
\def \P {{\cal P}}
\def \J {{\cal J}}
\def \F {{\cal F}}

\def \R {{\mathbf  R}}

\def\eps{\varepsilon}
\def\e{\varepsilon}

\usepackage{graphicx}
\graphicspath{{pictures/}}
\DeclareGraphicsExtensions{.pdf,.png,.jpg}
\begin{document}

\title{
Нерешенные задачи  о джойнингах, кратном перемешивании, спектре  и ранге   }
\author[V.\,V.~Ryzhikov]{В.\,В.~Рыжиков}
\address{Московский государственный университет}
\email{vryzh@mail.ru}

\date{2024}
\udk{517.987}

\maketitle

\begin{fulltext}

\begin{abstract}{ 

The note is devoted to multiple mixing, spectrum, rank and self-joinings of measure-preserving transformations. We recall famous open problems, discuss related questions and some known results.  A hypothetical example of an automorphism of the class $MSJ(2)\setminus MSJ(3)$ has Lebesgue spectrum, infinite rank and does not have multiple mixing. If its spectrum is simple, then we have a solution to the problems of Banach, Rokhlin and del Junco-Rudolph. The existence of such an example, of course, seems unlikely, but any facts confirming the impossibility of this amazing situation have not yet been discovered. They not found even under the  condition that its local rank is positive, which ensures the finite spectral multiplicity.

}
\end{abstract}

     
\markright{Нерешенные задачи о джойнингах}


\section{Спектр и   кратное перемешивание} 

\vspace{2mm}
\it  Предлагаемая работа возникла из попыток автора решить известную
спектральную проблему теории динамических систем: существуют ли
метрически различные динамические системы с одним и тем же непрерывным 
(в частности, лебеговским) спектром? С помощью введенных
для этой цели новых метрических инвариантов  автор  пытался найти 
среди эргодических автоморфизмов компактных коммутативных групп 
автоморфизмы различных метрических типов. Оказалось, однако, что у всех указанных автоморфизмов новые инварианты совершенно одинаковы.\rm

\vspace{2mm}
Выше мы процитировали первые строки  из  статьи Рохлина \cite{Ro}.
 Предложенный им инвариант  называются  
перемешиванием кратности $k\geq 1$. Напомним его определение. Автоморфизм $T$ вероятностного пространства $(X,\B,\mu)$ обладает свойством $k$-кратного перемешивания, если для всяких измеримых множеств $A,A_1,\dots, A_k$ при $m_1,\dots m_k\, \to\,+\infty$  выполняется
$$\mu(A\cap T^{m_1}A_1\cap T^{m_1+m_2}A_2\cap \dots\cap  T^{m_1+\dots+m_k}A_k)\ \to \ \mu(A)\mu(A_1)\dots \mu(A_k).$$
Обозначим это  свойство и класс соответствующих  автоморфизмов через $Mix(k)$.

\vspace{2mm}
{\bf  Проблема Рохлина о кратном перемешивании.} О непустоте классов 
$Mix(k)\setminus Mix(k+1)$ мы ничего не знаем. 

\vspace{2mm}
{\bf Q1. \it Существует  ли перемешивающий автоморфизм, не обладающий   перемешиванием некоторой кратности?} 

\vspace{2mm}
С.Каликов \cite{Ka} доказал перемешивание кратности 2 для  перемешивающих преобразований  ранга 1 (дальнейшие обобщения обсудим позже). Результат Каликова указывает на то, что двукратное  перемешивание совместимо со сколь угодно медленным обычным перемешиванием (мы, конечно,  умеем управлять скоростью перемешивания для преобразований ранга 1).  Общий  результат получил  Б. Ост:   при помощи джойнингов и тонких методов гармонического анализа им, доказано, что  автоморфизм из класса $Mix(k)\setminus Mix(k+1)$   обязан иметь абсолютно непрерывную компоненту в спектре \cite{Ho}. Можно удивляться и результату Оста и тому, как  промыслительно  Рохлин написал о лебеговском спектре в 1949 году.  

Хотя в этой стате мы интересуемся  задачами  о действии одного автоморфизма ($\Z$-действии), отметим,  что Ледрапье \cite{Le} дал элегантное решение  проблемы о кратном перемешивании в классе $\Z^2$-действий. Удивительно,  что  идея Рохлина использовать автоморфизмы компактных коммутативных групп оказалась  успешной для всех $\Z^n$-действий при $n>1$, причем примеры другой природы до сих пор  не обнаружены.

\vspace{2mm}
{\bf Проблема Банаха об однократном лебеговском спектре.} Проблема Рохлина  была упомянута  в списке из 10 задач в книге  П. Халмошем \cite{Ha}. Другая  нерешенная  задача  из этого списка известна как  проблема Банаха.

\vspace{2mm}
{\bf Q2.  \it Найдется ли автоморфизм пространства с мерой, обладающий однократным лебеговским спектром?}

\vspace{2mm}
 В этой  формулировке зашиты   две задачи. В книге Улама \cite{U} вопрос Банаха звучал так:  \it найдется ли сохраняющее меру Лебега $m$ на прямой $\R$ преобразование $T:\R\to\R$  и  функция $f\in L_2(m)$ такие, что векторы  $\{f\circ T^z\,:\, z\in\Z\}$ образуют полную ортонормированную систему?  \rm В случае пространства с вероятностной мерой вопрос аналогичен, но вместо $L_2$ рассматривается  пространство, ортогональное константам.  (Чтобы различать задачи,
последнюю можно  назвать  проблемой Банаха, а первую -- проблемой
Банаха-Улама.)

\vspace{2mm}
{\bf Q3. \it   Какими   спектрами могут обладать   эргодические автоморфизмы? Существует ли  автоморфизм с абсолютно непрерывным нелебеговским спектром?  Какие наборы спектральных кратностей могут быть у эргодического автоморфизма?}

\vspace{2mm}
Неизвестно, имеются ли эргодические автоморфизмы вероятностного пространства, у которых наборы спектральных кратностей имеют, например,  вид  $\{m, m+1\}$, $m>2$. 
Известно, что реализуется всякий  набор кратностей, содержащий  $1$, $2$ или $\infty$  (см. \cite{R23} и ссылки там на обзоры). Отметим, что все наборы спектральных кратностей реализуются эргодическими автоморфизмами пространства с сигма-конечной мерой  \cite{DR}.

\vspace{2mm}
{\bf Q4. \it Существует ли эргодические автоморфизмы с конечнократным абсолютно непрерывным спектром?}

\vspace{2mm}
Эту задачу  можно назвать  ослабленной проблемой Банаха.

\vspace{2mm}
В связи с проблемой Рохлина  был   упомянут результат Каликова: \it
если автоморфизм $T$ обладает перемешиванием и его локальный ранг $\beta (T)$ равен 1, то $T$ обладает перемешиванием кратности 2. \rm

\vspace{2mm}
Теорема Каликова  допускает усиление  на случай  $\beta (T)> \frac 1 2$ (\cite{93}).  Напомним определение инварианта $\beta (T)$.

\bf Локальный ранг автоморфизма. \rm
Говорят, что автоморфизм $T$ обладает \it локальным рангом \rm
$\beta >0$, если
для некоторой последовательности множеств вида
$U_j=\bigsqcup_{k\in Q_j}T^zB_j$,  $Q_j=\{0,1,\dots,h_j\}$, выполнены условия:
$\mu(U_j)\to \beta$, и для каждого $A\in\B$ пересечение  $U_j\cap A$
аппросимипуется объединениями некоторого набора этажей в башне: найдутся  $S_j\subset   \{0,1,\dots, h(j)-1\}$ такие, что 
$$\mu\left((U_j\cap A)\Delta\bigsqcup_{z\in S_j}T^zB_j\right)\ \to \ 0,  \ j\to\infty.$$ Через $\beta (T)$ обозначается максимальное из таких $\beta>0$, а если их нет, полагаем $\beta(T)=0$.

\vspace{3mm}
{\bf Одно из утверждений верно, хотя  не знаем, какое именно.} 
Следующая теорема  является следствием    результатов работ  Кинга, Тувено \cite{KT} и Глазнера, Оста, Рудольфа \cite{GHR}.

\vspace{2mm}
\bf  Теорема 1.1. \it  Верно по крайней мере  одно из следующих  утверждений:

-- существует автоморфизм   с конечнократным абсолютно непрерывным спектром (решние ослабленной проблемы Банаха);

-- всякий перемешивающий автоморфизм положительного локального ранга обладает кратным перемешиванием. \rm

\vspace{2mm}
{\bf Легкое кратное перемешивание.} Будем говорить, что перемешивающий автоморфизм $T$ обладает легким 
перемешиванием кратности 2, если для всяких последовательностей $m_i,n_i\to +\infty$ и множеств $A, B,C$ положительной меры  выполнено
$$ \limsup_{i\to\infty} \mu(A\cap T^{m_i}B\cap T^{m_i+n_i}C)> 0.$$
  \rm

\vspace{3mm}
\bf  Теорема 1.2. \it  Если перемешивающий автоморфизм $T$ положительного локального ранга не обладает перемешиванием кратности 2, то 
для некоторых  множеств $A,B$, $\mu(A), \mu(B)>0$ и последовательностей $m_i, n_i\to +\infty$  выполнено 
$$ A\cap T^{m_i}B\cap T^{m_i+n_i}B=\varnothing$$
(нет легкого  перемешивания кратности 2).\rm

\vspace{2mm}
Теорема 1.2 вытекает из результатов  \cite{97} о минимальных самоприсоединениях, которые мы изложим в \S 4.

\vspace{2mm}
\bf Q5. \it Влечет ли перемешивание легкое перемешивание кратности 2? \rm

\vspace{2mm}
{\bf Количество отклонений от кратного перемешивания.}
Пусть $\mu(A)>0$, положим 
$$D(T,A,N, a)= \{(m,n)\,:\, \mu(A\cap T^{m}A\cap T^{m+n}A)> 2 \mu(A)^3, 
\ aN<m,n<N\}.$$

\vspace{2mm}
\bf Теорема 1.3. (\cite{92}).  \it Для всякого перемешивающего автоморфизма $T$  выполнено $$\limsup_{N\to\infty} \frac {|D(T,A,N,a)|} N \ <\infty.$$ \rm

\vspace{3mm}
\bf  Теорема 1.4.  (\cite{95})\it  Если перемешивающий автоморфизм $T$ положительного локального ранга не обладает перемешиванием кратности 2, то количество отклонений от кратного перемешивания для него оптимально  в следующем смысле: для некоторых    $A$ и  $a>0$ выполнено 
$$\limsup_{N\to\infty} \frac {|D(T,A,N,a)|} N \ >0. \eqno (Max)$$ 
\rm

\vspace{2mm}
\bf Q6.  \it  Может   ли   свойство  $(Max)$ выполняться для перемешивающего автоморфизма?\rm

\vspace{2mm}
Предположение о том, что такие частые отклонения невозможны,
выглядит вполне разумным.

\vspace{2mm}
\bf О типичности кратного перемешивания. \rm Тихонов установил \cite{T12}, что типичные относительно  метрики Альперна-Тихонова  перемешивающие автоморфизмы обладают сингулярным спектром, следовательно,  в силу теоремы Оста они  обладают кратным перемешиванием.
Другое доказательство типичости кратного перемешивания, использующее теорему Каликова, получил Баштанов \cite{Ba}, доказав типичность ранга 1.  Иначе типичность ранга 1  может быть доказана с использованием перемешивающих орнстейновских преобразований (см. \cite{20}).  Изобилие задач  возникает в связи  с вопросом \it о типичности кратного перемешивания в пространствах перемешивающих групповых действий. \rm Отметим, что применение метода Бэра для групповых действий, оснащенных метриками типа  Альперна-Тихонова, дало полное решение проблемы Рохлина об однородном спектре в классе перемешивающих автоморфизмов \cite{Tih}.

Косые произведения (иначе говоря, расширения)  являются  источниками конструкций  в эргодической теории (в спектральной теории первые впечатляющие примеры были предложены Оселедецем без малого 60 лет тому назад). Идея решить проблему Рохлина при помощи расширений представляется весьма разумной, но все еще не реализована.   В \cite{2023} доказано, что типичное (в метрике Халмоша) расширение всякого перемешивающего автоморфизма также обладает  перемешиванием. Возникает   вопрос о кратном перемешивании.

\vspace{2mm}
\bf Q7.  \it Поднимается ли  свойство $Mix(k)$, $k>1$,  при типичном 
расширении? \rm

Возвращаясь к типичным перемешивающим автморфизмам, отметим, что 
они обладают свойством минимальных самоприсоединений (MSJ),  введеным Рудольфом в \cite{Ru}. Это свойство воспринималось  как экзотика:
оно влечет за собой  явное описание  минимальной структуры факторов и централизаторов всех декартовых степеней автоморфизма.  Оказалось, что  экзотика  типична, но не в метрике Халмоша, а в  метрике Альперна-Тихонова.

\section{Джойнинги} 
Cвойства действий, формулируемые в терминах  джойнингов,   оказались полезными    при  исследовании  кратного перемешивания (и многих других свойств динамических систем, см., например, \cite {Th},\cite{Gl}). Приведем необходимые определения. 

Присоединением (или джойнингом) набора действий $\Psi_1,\dots,\Psi_n$ называется мера на $X^n =X_1\times \dots\times  X_n$  ($X_i=X$), проекции которой на ребра  куба $X^n$ равны $\mu$, причем мера инвариантна относительно диагонального действия произведения $\Psi_1\times \dots\times  \Psi_n$. Если $\Psi_1,\dots,\Psi_n$ суть копии одного  действия, такой джойнинг называется самоприсоединением (self-joining).

Говорим, что действие $\Psi$  принадлежит классу $S(m, n)$, $n>m>1$  (или обладает свойством $S(m,n)$), если всякое cамоприсоединение  порядка $n>2$ такое, что все проекции на  $m$-мерные грани  куба $X^n$ равны $\mu^m$,  является тривиальным,  т.е. совпадает с мерой $\mu^n$,
произведением $n$ копий меры $\mu$.  Далее используем обозначение   $S_n=S(n-1,n)$, $n>2$. Отметим, что положительная энтропия и наличие дискретного спектра несовместимы со свойствами $S_n$.  

{\bf  Кратное перемешивание, джойнинги и спектр.} Ряд результатов о кратном перемешивании был получен без   кропотливого изучения  кратных персечений  как следствие некоторых  свойств джойнингов.

\vspace{4mm}
\bf Лемма 2.1. \it Свойства $S_n$  связаны с $Mix(k)$ следующим образом: если перемешивающий автоморфизм обладает перемешиванием кратности $k$ и свойством $S_{k+2}$, то он обладает кратным перемешиванием  порядка $k+1$. \rm

\vspace{2mm}
Поясним, как устанавливается указанная связь. Пусть для всяких измеримых $A,B,C$ выполняется
$$\mu(A\cap T^{m_i} B\cap T^{m_i+n_i}C)\to \nu(A\times B\times C), $$
Свойство $S_3$  автоморфизма $T$ влечет за собой
$$\nu(A\times B\times C) =\mu(A)\mu(B)\mu(C).$$
Действвительно,  $\nu$ --  нормированная мера на полукольце цилиндров 
вида $A\times B\times C$.  
Проверяется, что  
$$\nu(A\times B\times C) = \nu(TA\times TB\times TC),$$ 
$$\nu(A\times B\times X)=\mu(A)\mu(B), \ \nu(X\times B\times C) = \mu(B)\mu(C),$$
$$ \nu(A\times X\times C) = \mu(A)\mu(C).$$
Такие меры называют самоприсоединениями  с попарной независимостью.

В связи с этой леммой естественно задаться вопросом "бывают ли перемешивающие преобразования со свойством $S_n$? Ответ содержит  работа Рудольфа \cite{Ru}.   В этой статье Рудольф параллельно устанавливал свойство MSJ (которое по определению сильнее всех свойств $S_n$) и доказывал кратное перемешивание, не обращая внимания на связь этих инвариантов (автор заметил эту связь в 1986 году и увлекся джойнингами).  В  статье  \cite{JR} дель Джунко и Рудольф   ввели инвариант   PID, означавший  выполнению всех свойств $S_n$, показали его устойчивость относительно прямых произведений и компактных расширений, упомянули  связь инварианта   с кратным перемешиванием. Результаты Каликова и Кинга стимулировали новые наблюдения о  совпадении свойств  MSJ(2) и MSJ для потоков и доказательство    свойств $S_n$ и $Mix(k)$ для перемешивающих преобразований конечного ранга \cite{91}.

  Глубокий  результат  о свойствах $S_n$ для автоморфизов был получен Остом.

\vspace{3mm}
\bf Теорема 2.2. (\cite{Ho}). \it Автоморфизм с непрерывным сингулярным спектром обладает всеми свойствами $S_n$. \rm

\vspace{2mm}

\vspace{2mm}
{\bf  Взаимосвязи свойств $S_n$. Четное и нечетные свойства.}
Кинг \cite{Ki}  обнаружил, что  \it одновременное выполнение свойств $S_3$ и $S_4$ влечет за собой выполнение всех свойств $S_p$.  \rm
 
В  \cite{97S}  доказано, что  \it четные свойства $S_{2m}$ для всех $m>1$ эквивалентны между собой и  вынуждают  все нечетные свойства  $S_{2n-1}$, $n>1$.  
\rm 

\vspace{2mm}
\bf Q8. \rm Как связаны между собой нечетные свойства  $S_{2n-1}$?

\vspace{2mm}
Известны  групповые действия класса   $S_{2m-1}\setminus S_4$.  Пример: рассмотрим компактную коммутативную группу $X=\Z_2^{\Z}$ с мерой Хаара $\mu$,  определим на $(X,\mu)$  действие $\Psi$, порожденное всеми сдвигами на группе $X$ и всеми автоморфизмами группы $X$. Действие $\Psi$  и многие его  поддействия обладают всеми нечетными свойствами  $S_{2m-1}$, но не обладают четным свойством  $S_{4}$ (см. ссылки в \cite{97S}). 

\vspace{2mm}
{\bf Свойство $\bf S_3$ и наследственная независимость факторов.}  Действие $T$ по определению обладает НН свойством (НН от \it Наследственная Независимость \rm), если для любой  эргодической  динамической  системы, порожденной тремя факторами
${\cal F}, {\cal F}',{\cal F}''$, каждый из которых изоморфен системе $T$, выполнено
$$
    {\cal F}\bot{\cal F}' \ \ \& \ \  {\cal F}\bot{\cal F}'' \ \
\Rightarrow \ \ \ {\cal F}\ \bot\ ({\cal F}'\bigvee {\cal F}''),
$$
\rm 
где $\bot$ обозначает независимость факторов.

Заметим, что НН  влечет в силу определения свойство $S_3$.
Оказывается \cite{97}, что НН влечет $S_4$. Мы изложим в \S 3 доказательство этого утверждения. 

{\bf  Связь с энтропией. Проблема дель Джунко-Рудольфа \cite{JR}.}
О понятии энтропии автоморфизма вероятностного пространства см. \cite{KSF}.

\vspace{2mm}
\bf Q9.   \it Существует ли слабо перемешивающий автоморфизм с нулевой энтропией, не обладающий свойством $S_n$? \rm

\vspace{2mm}
{\bf Независимый фактор. }  Не решен  следующий весьма  частный случай этой задачи. 

\vspace{2mm}
\bf Q10. \it
Пусть алгебра $\F$ (фактор) является  инвариантной относительно  эргодического произведения $T\otimes T$ и  независима от координатных алгебр. Верно ли, что в этом случае энтропия автоморфизма $T$ положительна?\rm

\vspace{2mm}
\begin{center}
\includegraphics{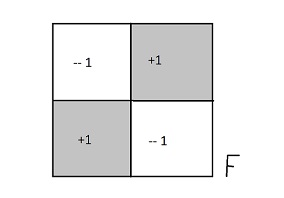}
\end{center}

\bf Теорема 2.3. (\cite{92}) \it Если фактор $\F$ порождается множеством вида $(A\times A) \cup (\bar A\times\bar A)$ , 
$\mu(A)=1/2$, $\bar A$ -- дополнение к $A$, то такой фактор  бернуллиевский. \rm

\vspace{2mm}
Доказательство.  Обозначим  $f=2\chi_A -1$, $F_n= (T^nf\otimes T^nf)$ (на  рисунке $F=F_0$).
Покажем, что $T^nf$ -- последовательность независимых случайных величин.
Замечаем, что  $F_n \in L_2(\F)$ (измеримы относительно фактора $\F$. 
 Пусть  $n_1,\dots, n_k$ различны, положим 
$$c:=\int T^{n_1}f\dots T^{n_k}f\, d\mu,$$ 
тогда
$$c^2=\int F_{n_1}\dots F_{n_k}\, d(\mu\times \mu).$$ 
Имеем
$$c=\int \un\otimes T^{n_1}f\dots T^{n_k}f\,d(\mu\times \mu)  = $$
$$=\int (F_{n_1}\dots F_{n_k})(T^{n_1}f \dots \bar T^{n_k}f\otimes \un) \, d(\mu\times \mu)  =c^3,$$ 
последнее равенство вытекает из независимости фактора $\F$  от координатного фактора.
Получаем $c=0,1,-1$.  Если $T$ слабо перемешивает, возможен только случай $c=0$, как в случае $c=1,-1$  все $F_n$ лежат в  конечномерном пространстве.
Например, если $$-1=F_0\, F_3\, F_{10}F_{11},$$  тогда 
$$F_{11}=-F_0\, F_3\, F_{10}, \ \ 
F_{12}=-F_1\, F_4\, F_{11}= F_1\, F_4\, F_0\, F_3\, F_{10}$$ 
$$F_{13}= -F_2\, F_{5}\, F_{12}, \ \  
F_{14} = F_3\, F_6\, F_{13}, \dots .$$
Таким образом,   всякая функция $F_n$ при $n>10$ измерима относительно конечной алгебры множеств, порожденной функциями $F_1,\dots, F_{10}$, так как выражается через их произведения. Общий случай аналогичен.

Равенство нулю указанных интегралов означает независимость случайных величин  $F_n$ и тем самым  независимость функций $T^n f$.
В нашем  случае независимый  фактор является схемой Бернулли типа $(\frac 1 2, \frac 1 2)$ и   мы получили  эффект положительности энтропии автоморфизма $T$.  

Всякое эргодическое преобразование с положительной энтропией имеет бернуллиевский фактор \cite{KSF} и по этой причине  не обладает свойствами $S_n$. Большие  степени  декартового квадрата автоморфизма всегда обладают независимым фактором, который мы только что описали.

\vspace{2mm} 
{\bf Гауссовские автоморфизмы. }
 Для любой непрерывной борелевской меры $\sigma$ на единичной окружности в комплексной плоскости найдется эргодический гауссовский автоморфизм (см. \cite{KSF}), спектр которого содержит $\sigma$   в качестве компоненты.
Для некоторых мер такой автоморфизм в определенном  смысле единственный.   Если же мера $\sigma$  сингулярна, то соответствующий гауссовский автоморфизм имеет нулевую энтропию. В связи с этим уместно  обратить внимание читателя на  частный случай проблемы дель Джунко-Рудольфа.

\vspace{2mm}
\bf Q11. \it Обладают ли свойствами $S_n$ эргодические гауссовские автоморфизмы с нулевой энтропией ?\rm

\vspace{2mm}
{\bf Минимальные самоприсоединения.}  Рудольф ввел в эргодическю теорию свойство автоморфизма,   контрастирующее со свойствами гауссовских систем.  Речь идет о преобразованиях с минимальными самоприсоединениями (их класс обозначается через $MSJ$ как абревиатура от  minimal self-joinings).   Определим при $n>1$  класс $MSJ(n)$ так:  пусть централизатор автоморфизма $T$ тривиален, а для всякой эргодической системы и любого  набора из $n$  факторов, изоморфных автоморфизму $T$, часть из них образуют независимую систему, а остальные  совпадают с упомянутыми. Класс $MSJ$  определен как пересечение всех $MSJ(n)$, $n>1$.  В  \cite{GHR}, в частности, доказано, что   $MSJ(3)=MSJ$.

\vspace{2mm}
{\bf Можно ли решить одновременно проблемы Банаха, Рохлина и  дель Джунко-Рудольфа? } Гипотетический пример автоморфизма класса  $MSJ(2)\setminus  MSJ(3)$ обладает лебеговским спектром \cite{GHR},  бесконечным рангом  \cite{93} и  не обладает кратным перемешиванием \cite{97}. А если  его спектр простой, это дает решение проблем Банаха, Рохлина и  дель Джунко-Рудольфа.  Существование такого  примера, конечно,   представляется маловероятным, но факты, подтверждающие невозможность  удивительной ситуации, пока не обнаружены. Они не найдены за 30 лет даже при дополнительном условии положительности  локального ранга автоморфизма, которая гарантирует конечную кратность спектра. 

\vspace{2mm}
\bf Q12.  \it  Cовпадают ли классы MSJ(2) и MSJ(3)? \rm

Для неперемешивающих преобразований совпадение классов 
$MSJ(2)$ и $MSJ(3)$ вытекает, например, из результатов \cite{GHR}.
Ответ положителен для перемешивающих автоморфизмов конечного ранга и автоморфизмов
с локального ранга, превосходящего $\frac 1 2$. Такие автоморфизмы   обладают так называемым  D-свойством, которое 
влечет свойство $S_4$, следовательно, свойство MSJ(4). 
 
{\bf D-свойство.} Последовательность разбиений множеств $U_j\subset X$
вида
$$ \xi_j=\{ E_j, TE_j,\dots T^{h_j}E_j\}$$
назовем   аппроксимирующей,
если, дополняя разбиение $ \xi_j$ некоторым разбиением  дополнения
$X\setminus U_j$,
получим последовательность разбиений всего фазового пространства $X$,
которая стремится к разбиению на точки.

Будем говорить, что автоморфизм $T$ обладает D-свойством, если
найдутся последовательности аппроксимирующих башен
$ (U_j,\xi_j)$, $(U'_j,\xi'_j)$, $(U''_j,\xi''_j)$, 
где
$$ \xi_j=\{ E_j, TE_j,\dots T^{h_j}E_j\},\
 \xi'_j=\{ E'_j, TE'_j,\dots T^{h_j}E'_j, \
\xi''_j=\{ E''_j, TE''_j,\dots T^{h_j}E''_j\},$$
причем  для некоторой последовательности $m_j$, $m_j>h_j$,
выполняются следующие условия:
$$\lim_j \mu(U_j) =a>0, \ \  \mu(E_j)=\mu(E'_j)=\mu(E''_j), $$
$$ E'_j= T^{m_j}E_j,\ \  \mu(T^{m_j}U'_j\Delta U''_j)\to 0, $$
$$ \max_{m>h_j}\mu(T^{m}E'_j\ |\ E''_j)\ \to\ 0.$$

\vspace{2mm}
{\bf Теорема 2.4. (\cite{93}) } \it
 Перемешивающий автоморфизм конечного ранга и перемешивающий автоморфизм локального ранга $\beta>\frac 1 2$ 
обладают D-свойством.
  Перемешивающий автоморфизм c   D-свойством принадлежит классам $S_4$ и $Mix(k)$,  $k>1$. \rm

\section{Задача  $\bf S_3 \Rightarrow S_4$ и наследственная независимость
факторов}
В этом параграфе, следуя \cite{97} мы изложим обобщение результата \cite{GHR} о том, что 
свойство простоты порядка 3 влечет за собой простоту порядка 4.

Действие $T$ по определению обладает свойством НН, если для любой
эргодической  динамической  системы, порожденной тремя факторами
${\cal F}, {\cal F}',{\cal F}''$, каждый из которых изоморфен системе $T$, выполнено
$$
    {\cal F}\bot{\cal F}' \ \ \& \ \  {\cal F}\bot{\cal F}'' \ \
\Rightarrow \ \ \ {\cal F}\ \bot\ ({\cal F}'\bigvee {\cal F}''),
$$
\rm
(где $\bot$ обозначает независимость факторов).

Переформулируем  это определение.
Действие $T$ обладает свойством НН, если \it для любого эргодического
джойнинга  $\eta$ тройки  $(T,T,T)$  выполнено
$$
    \pi_{12}\eta =\mu\otimes\mu \ \ \& \ \ \pi_{13}\eta =\mu\otimes\mu
\ \ \Rightarrow \ \ \eta = \mu\otimes \pi_{23}\eta,
$$
где  $\pi_{ij}\eta$ --  проекция
 меры  $\eta$
на грань  $X_i\times X_j$ куба $X_1\times X_2\times X_3.$
\rm

\vspace{2mm}
{\bf Теорема 3.1} \it Автоморфизм с 
свойством НН  принадлежит классу  $S_4$.\rm

\vspace{2mm}
 Пусть $\nu$ -- попарно независимый джойнинг набора $(T,T,S)$,
где $S$ -- некоторый эргодический
автоморфизм, а  $T$ является  автоморфизмом со свойством НН.

Для доказательства  теоремы  мы воспользуемся индуцированными
джойнингами: мера $\nu$ индуцирует последовательность
 джойнингов  $\eta_m$. Дадим их определение.

Пусть семейство марковских операторов $\{P_x\}$ отвечает мере
$\nu$ :
для всех $A,B,C\in\B$
$$
  \int \! \chi_A(x)\langle P_x\chi_B ,\chi_C \rangle d\mu(x) \
= \ \nu (A\times B\times C),
$$
где $\chi_A $ -- индикатор  множества $A$,
$\langle\ , \ \rangle$ обозначает скалярное произведение
в пространстве $L_2(X,\mu)$.
Из инвариантности  меры  $\nu$ относительно $T\times T\times T$
вытекает тождество
$$S^{-1}P_{T^{-1}(x)}T \equiv P_x .$$

 Рассмотрим джойнинги $\eta_m,$ $ m\in Z,$, индуцированные мерой $\nu$,
 заданные  равенствами
$$
    \eta_m (A\times B\times C) = \int \! \chi_A(x)
\langle P_x\chi_B ,P_{T^{m}x}\chi_C \rangle d\mu(x).
$$
 Мера $\eta_m$  инвариантна относительно
$T\times T\times T$.
\medskip

{\bf Лемма  3.2. }\it
Если $T$ является HI-системой, то равенство
$$
P_{T^{m}x}^\ast P_x = \int_X P_{T^{m}x}^\ast  P_x d\mu(x)
$$
выполнено для почти всех  $x$.
\rm
\medskip

{ Доказательство}.
Фиксируем $m$, пусть $\eta =\eta_m$. Из определения
$\eta$ имеем
$$
\eta (B\times C\times X) = \int \chi_B d\mu\int \chi_C d\mu
=\mu(B)\mu(C)= \eta (B\times X\times C),
$$
следовательно,
$$
    \pi_{12}\eta =\mu\otimes\mu, \  \ \pi_{13}\eta =\mu\otimes\mu.
$$
Так как $T$  есть HI-система, мы получаем
$$
\eta = \mu\otimes\ \pi_{23}\eta,
$$
что влечет за собой выполнеие  равенств
$$\eta_{x}=\pi_{23}\eta$$
для почти всех условных  мер $\eta_{x}$ (соответствующих
 операторам $P_{T^{m}x}^\ast P_x$). Поэтому
для почти всех $x\in X$ выполнено
$$P_{T^{m}x}^\ast P_x = \int_X P_{T^{m}x}^\ast  P_x d\mu(x).$$
\\  

 Доказательство теоремы 3.1.     Последовательность $m_i\to\infty$ будем называть  перемешивающей (дразумевая перемешивающей
последовательность степеней $T^{m_i}$), если
$$
\forall A,B\in \B \ \ \ \mu(A\cap T^{m_i}B) \to  \mu(A)\mu(B).
$$
Существование  перемешивающей последовательности
вытекает из свойства  слабого  перемешивания. Ниже $\Theta$ обозначает ортопроекцию на  константы в $L_2(\mu)$.
\medskip

\bf Лемма  3.3. \it
Для перемешивающей последовательности $\{m_i\}$ имеет  место
следующая слабая операторная сходимость:
$$
\int_X P_{T^{m_i}(x)}^\ast P_x  d\mu(x) \to \Theta.
$$
\rm
\medskip

{ Доказательство.}
Для заданных функций $f,f'\in L_2(\mu)$,
$\ 0\leq f(x),\, f'(x)\leq 1,$ и $\eps >0$
найдется  набор
марковских операторов $\{P_1,\dots , P_N\}$  такой, что
фазовое пространство $X$ можно представить в виде
$X=B\cup\cup_{k=1}^N A_k$, где $\mu(B)<\eps$ и
$
  \forall k \ \forall x,x'\in A_k $
$$\|P_xf-P_kf\|<\eps,
  \ \ \|P_xf'-P_kf'\|<\eps.
$$
Тогда с учетом свойство перемешивания
$$\mu(T^{-m_i}A_k\cap A_l)\to \mu(A_k)\mu(A_l)$$
и равенства $$\int_x P_x d\mu(x)=\Theta$$
  получаем:
$$
\int_X P_{T^{m_i}(x)}f P_x f' d\mu(x) \ \leq  \
\sum_{k,l}\int_{T^{-m_i}A_k\cap A_l} P_{T^{m_i}(x)}f P_x f' d\mu(x)
\   +   \ 2\eps  \ \leq
$$
$$
\leq \sum_{k,l}\int_{T^{-m_i}A_k\cap A_l} P_kf P_l f' d\mu(x)
 \  +    \ 4\eps  \ \leq  $$
$$
\leq \sum_l\int_{A_l}(\sum_k \mu(A_k)P_kf) P_l f' d\mu(x)
\   +   \ 5\eps  \ \leq
$$
$$
\leq \sum_l\int_{A_l}(\int fd\mu) P_l f' d\mu(x)
\   +    \ 6\eps  \ \leq
\int_X fd\mu \int_X f'd\mu
\   +  \   7\eps.
$$
\medskip

\bf Лемма   3.4. \it Для почти всех  $y$   найдется  подпоследовательность
$m_i(y)\to\infty$ такая, что имеет место слабая сходимость
$$
 P_{T^{m_i(y)}(y)}^\ast P_y  \to  P_y^\ast P_y.
$$
{ Доказательство.}     \rm
Для любых $\eps >0$ и $f\in L_2(\mu)$
найдутся множества
$A_1, \, A_2,\,\dots, A_p,\dots $ положительной меры
такие, что
$$
\forall p \ \ \  \sup_{x,y\in A_p}  \|P_xf-P_yf\| \leq \eps.
$$
Так как  $\{ m_i\}$ является  перемешивающей последовательностью,
для почти всех $y$ найдется  подпоследовательность
 $\{m_{i'}\}$  (зависящая от $y$)
такая, что $T^{m_{i'}}(y)\in A_p$. Тогда выполнено
$$
         \|P_{T^{m_{i'}}(y)}f-P_yf\| \ \leq \eps.
$$
Используя обычную диагональную процедуру и сепарабельность
пространства $L_2(\mu)$), получим
$$ P_{T^{m_{i''}}(y)}^\ast P_y \to  P_y^\ast P_y.$$
\\ 
Завершим доказательство теоремы 3.1.
В силу  леммы 3.3 выполнено
$$
\int_X P_{T^{m_i(y)}(x)}^\ast P_x  d\mu(x) \to \Theta ,
$$
а из лемм 3.2. и 3.4  вытекает, что для почти всех
$y$ выполнено $\ P_y^\ast P_y =\Theta$, $P_y =\Theta$. Следовательно,
$\nu=\mu\otimes\mu\otimes\mu$.


 \section{Задача  $\bf MSJ(2) \Rightarrow MSJ(3)$ и  легкое кратное перемешивание}
В этом параграфе, следуя \cite{97}, мы излагаем  решение задачи 
$MSJ(2)=MSJ(3)?$ при  дополнительном условии типа  перемешивания кратности 2.

Автоморфизм $T$ пространства $(X,\B,\mu)$
 обладает свойством минимальных самоприсоединений
порядка
$n$ ($T\in MSJ(n)$), если
любой эргодический  джойнинг  $n$ копий $T$,
исключая меру $\mu^{\otimes n}=\mu_{(1)}\otimes\dots\otimes\mu_{(n)}$,
обладает следующим свойством:
одна из его проекций  на двумерную грань  в
$X\times\dots\times X$ являеся мерой  $\Delta_{T^i}$ (сдвиг диагональной
меры).  Неформально говоря,  такой автоморфизм $T$ имеет только
очевидные джойнинги.

Для  $\Z$-действий мы  покажем, что  в классе
$MSJ(2)$ проблема Рохлина эквивалентна  открытому
вопросу  терии джойнингов:
\it совпадает ли класс $MSJ(2)$ с классом  $MSJ(3)$? \rm

Хотя имеются некоммутативные контрпримеры, т.е. для групповых действий
возможно, как мы показали в главе 1, несовпадение классов
$MSJ(2)$ и $MSJ(3)$ (и даже $MSJ(3)\neq MSJ(4)$),
для  $\Z$-действий, как отмечалось, задача не решена.
\medskip

{\bf Теорема 4.1.} \it  Если $T\in MSJ(2)$ и $T$ перемешивает с кратностью 2, то
автоморфизм $T$ обладает минимальными самоприсоединениями всех порядков и
кратным перемешиванием всех порядков.
\rm
\medskip

Этам теорема является непосредственным следствием более общего утверждения.
\medskip

{\bf Теорема 4.2.} \it   Пусть перемешивающий
 автоморфизм $T\in MSJ(2)$ обладает свойством
кратного возвращения:
  для любого множества $A$ положительной меры  и любой
 последовательности  $k(m)\to \infty$,
$|k(m)-m|\to \infty$  для всех больших $m$ выполнено
условие
$$
 \mu(T^{-k(m)}A\cap T^{-m}A\cap A) > 0.
$$
Тогда $T\in MSJ(3)$ и, следовательно, обладает свойством кратного
перемешивания всех порядков.
\rm
\medskip

ЗАМЕЧАНИЕ.  Сформулированная теорема, в частности, утверждает
следующее:
для перемешивающего
автоморфизма $T\in MSJ(2)$ свойство
$$\lim_{m\to\infty} \mu(T^{-k(m)}A\cap T^{-m}A\cap A) > 0$$
для любого $A$, $\mu(A) > 0$
влечет за собой
$$\mu(T^{-k(m)}A\cap T^{-m}A\cap A) \to \mu(A)^3.$$
\medskip

Сейчас мы сформулируем техническое утверждение,
играющее ключевую роль в доказательстве.

\it  {\bf Утверждение  4.3.}
Пусть автоморфизм $T$ принадлежит классу $ MSJ(2)\setminus MSJ(3)$,
тогда найдется число $a>0$,  множество $M\subset \N$
положительной  плотности,   семейство марковских операторов
 $\{\J_x\},\, \J_x:L_2(\mu)\to L_2(\mu)$, отвечающих
некоторому нетривиальному эргодическому
джойнингу со свойством парной независимости,
и  последовательность $k(m)$ такая, что
  $k(m)\to \infty$,   $|k(m)-m|\to \infty$ и
для любых множеств $A',B$ положительной
меры  для некоторого $A\subset A'$, $\mu(A)>0$
неравенство
$$
   \langle \J_{T^{k(m)}(x)} \chi_B |\chi_B \rangle
     \ >\ a\mu(B)
$$
выполнено для всех   $x\in {A\cap T^{-m}A}$ при  $m\in M$.
\rm
\medskip


Утверждение будет доказано позже, пока мы выведем из него
теорему.

 Доказательство теоремы 4.2.
Пусть $T\in MSJ(3)$ не выполняется. Ввиду утверждения 4.3 имеем:
для любого фиксированного $\e >0$ и   множества $B\in \B$  пространство
$X$ представляется как объединение
некоторых дизъюнктных множеств $A_1,A_2, \dots ,$ таких, что   выполнено
$$
\forall j\   \forall x,x'\in A_j \ \ \  \|\J_x\chi_B -\J_x'\chi_B\|<\e,
$$
причем для всех точек $x\in A_j\cap T^{-m}A_j$ имеет место неравенство
$$
\langle \J_{T^{k(m)}(x)}\chi_B |\chi_B \rangle
     \ >\ \frac{a}{2} \mu(B).
$$
Кратное возвращение обеспечивает следующее:
найдется $x\in A_j\cap T^{-m}A_j$ такая, что ${T^{k(m)}(x)}\in A_j$.
Поэтому выполнено
$$
\forall x'\in A_j \ \ \
\langle  \J_{x'}\chi_B \,|\,\chi_B \rangle\ \ > \ \frac{a}{2} \mu(B) - \eps.
$$
Таким образом, при $\frac{a}{10}\mu(B)>\eps>0$  и
$\mu(B)< \frac{a}{2}$ мы получаем  противоречие:
$$
 \mu(B)^2 = \int_X \langle  \J_{x'}\chi_B \,|\,\chi_B \rangle d\mu(x') =
 \langle  \Theta\chi_B\, |\,\chi_B \rangle
     \ \geq\ \frac{a}{2}\mu(B).
$$
Следовательно, предположение $T\in MSJ(2)\setminus MSJ(3)$ неверно.

Теперь приступим к доказательству утверждения 4.3.

Основная  идея состоит в следующем.
Индуцированный  джойнинг $\eta_m$
(индуцированный некоторым эргодическим джойнингом $\nu$ с попарной
независимостью) имеет следующее представление:
     $$
       \eta_m = \frac{1}{p}( \nu_{i(m,1)}+\nu_{i(m,2)}+\dots +\nu_{i(m,p)}),
     $$
где $\nu_{i(m,j)}$ --  эргодические джойнинги класса $M(2,3)$
(конечно, при  $m\neq 0$).
Можно доказать, что одна из компонент, скажем, $\nu_{i(m,1)}$,
стремится к $\Theta$. Это влечет за собой тривиальность исходного
джойнинга  $\nu$.

\newpage
{\bf Лемма   4.4. }\it    Пусть $\nu$ -- эргодический джойнинг набора
$(T,T,T)$ и выполнены условия: $\nu\in M(2,3)$,
$\nu\neq \mu\otimes\mu\otimes\mu$
 и $T\in MSJ(2)$.
Пусть $\{\P_x\}$ -- марковский  оператор, отвечающий   мере $\nu_x$,
где $\{\nu_x\}$ ($x\in X_{(1)}$) -- семейство  условных мер
на $X_{(2)}\times X_{(3)}$, отвечающих
джойнингу $\nu$.

Тогда найдутся целые числа $p,q\geq 1$ такие, что

(i) для почти всех $x$ выполнено
$$\P_x^\ast\P_x \ \geq \ \frac{1}{q}I; $$

(ii)
для всех $m\in\N, \ m\neq 0$ имеет  место тождество
     $$
       \P_{T^m(x)}^\ast \P_x  = \frac{1}{p}
          ( \J_x^{i(m,1)}+\J_x^{i(m,2)}+\dots +\J_x^{i(m,p)});
     $$

(iii)
Пусть $B$ --  множество  положительной  меры.
 Для любого множества $A'$, $\mu(A')>0$  найдется множество $A\subset A'$
положительной меры  такое, что для всех $x\in A\cap T^{-m}A$
выполнено

$$
  \langle \P_{T^m(x)}^\ast \P_x \chi_B |\chi_B \rangle
     \ >\ \frac{1}{2q} \mu(B);
$$

(iv)
   для любого $m$ найдется  число $r(m), \ 1\leq r(m)\leq p$, такое, что
для всех $x\in A\cap T^{-m}A$
$$
  \langle \J_x^{i(m,r(m))}\chi_B |\chi_B \rangle
     \ >\ \frac{1}{2q} \mu(B);
$$

(v)     для некоторого множества $M\subset \N$  положительной плотности
для  эквивариантного семейства
        $\{\J_x\}$, соответствующего некоторому джойнингу класса $M(2,3)$,
        выполнено
$$
 \forall m\in M \ \ \
   \J_x^{i(m,r(m))} =\J_{T^{k(m)}(x)}.
$$
\rm
\medskip

{ Доказательство} пункта (i).

 Равенство
  $$\int \P_x^\ast\P_x d\mu(x)= \Theta$$
влечет за собой
$$\P_x^\ast\P_x \equiv \Theta,  \ \ \P_x \equiv \Theta, \ \
\nu=\mu\otimes\mu\otimes\mu. $$
Так как $T\in MSJ(2)$, то оператор   $\int \P_x^\ast\P_x d\mu(x)$ является
выпуклой суммой оператора $\Theta$ и степеней $T^i$.

Так как $\nu\neq\mu\otimes\mu\otimes\mu$,
 для некоторого целого $m$ и числа $a>0$ имеем
$$\int \P_x^\ast\P_x d\mu(x) = aT^i +\dots .$$
Заметим, что случай $i\neq 0$ невозможен. Действительно,
из $\int \P_x^\ast\P_x d\mu(x) \geq aT^i$ вытекает, что
для почти всех $x$ операторы $\P_x$ и $\P_xT^i$
имеют ``общую часть'', т.е.  мера
$\nu$ и мера $(Id\times Id\times T^i)\nu$ имеют
общую компоненту. Но эти меры эргодичны относительно
$T\times T\times T$, следовательно, они совпадают.
Таким образом, мы получили равенство
$$\nu=(Id\times Id\times T^i)\nu. $$
Пусть $i\neq 0$.  Так как $\nu$ принадлежит $M(2,3)$, а
 $T^i$ -- эргодический автоморфизм (автоморфизм класса $MSJ(2)$ обязан быть
слабо перемешивающим), мы получаем $\nu =\mu\otimes\mu\otimes\mu$
(см. принцип дополнительно симметрии).
Таким образом, возможен только  случай $i=0$.
\medskip

{ Доказательство} пункта (ii). Имеем
$$\P_x^\ast\P_x \geq \frac{1}{q}I;
\ \ \P_x\P_x^\ast  \geq \frac{1}{r}I. $$
Тогда для $ {\cal H}_{mx}= \P_{T^m(x)}^\ast \P_x$ мы получаем  неравенства
$$ {\cal H}_{mx}^\ast {\cal H}_{mx} \geq \frac{1}{qr}I,
\ \  {\cal H}_{mx} {\cal H}_{mx}^\ast \geq \frac{1}{qr}I . $$

Теперь  для фиксированного $m$ рассмотрим  систему
$(T\times T\times T, \eta_m)$, где
$\eta_m$ --  джойнинг,  отвечающий семейству $\{ {\cal H}_{x}\}$.
Эта система для некоторого $s\leq qr$ является  $\Z_{s}$-расширением
системы    $(T\times T, \mu\otimes \mu)$.

Так как  последняя  эргодична,  число эргодических компонент
системы $(T\times T\times T, \eta_m)$ не превосходит числа $qr$.
Все компоненты являются джойнингами класса $M(2,3)$.
 Иначе мы получим
$$\nu=(T^m\times Id\times T^i)\nu,  \ m\neq 0,$$
но это влечет за собой равенство $\nu=\mu \otimes\mu \otimes\mu$,
так как  джойнинг $\nu$ принадлежит классу $M(2,3)$, а автоморфизм
$T^m\times T^i$
эргодичен при $i\neq 0$.
\medskip

{ Доказательство} пункта (iii).

    Зафиксируем некоторые множества $A'$ и $B$  положительной меры, для
которых выполнено условие
$P_x\chi_B\neq const$ при
$x\in A'$. Обозначим $ \hat{B}=\chi_B-\Theta\chi_B$.
 Рассмотрим  множество $A'$   положительной меры  такое, что
для некоторого $c,\ 0< c <\frac{1}{q}$  неравенство
$\|P_x \hat{B}\|>c$ выполнено  для всех $x\in A'$.  Для  $\eps>0$
выберем  $f\in L_2(\mu)$
и   множество $A\subset A'$ положительной меры
такое, что  $\|P_x \hat{B}- f\| < 0.1c\mu(B)$ для всех $x\in A$.
Существование  такой   функции $f$ следует из
 сепарабельности  пространства $L_2(\mu)$.
Получаем
$$
 \forall x,x'\in A\ \ \ \|\P_x\hat{B}- \P_{x'} \hat{B}\| < 0.2c\mu(B),
$$
следовательно,
$$
 \forall x,x'\in A\ \ \ \
\|\P_{x'}^\ast\P_x \hat{B}- \P_{x}^\ast\P_x \hat{B}\| < 0.2c\mu(B).
$$
Поскольку выполнено
$$\P_x^\ast\P_x \chi_B\geq \frac{1}{q}\chi_B, $$
при $c <\frac{1}{q}$ мы получаем
$         \forall\  x\in A\cap T^{-m}A$
$$
  \langle \P_{T^m(x)}^\ast \P_x \chi_B |\chi_B \rangle
   \ >\ \frac{1}{q} \mu(B) - 0.2c\mu(B)  \ >\ \frac{1}{2q} \mu(B).
$$
\medskip

{ Доказательство} пункта (iv).

Из пункта (ii) вытекает представление
     $$
       \P_{T^m(x)}^\ast \P_x  = \frac{1}{p}
          ( \J_x^{i(m,1)}+\J_x^{i(m,2)}+\dots +\J_x^{i(m,p)}).\eqno (4.1)
     $$
С точностью до перестановки членов в приведенной выше сумме
для всех $x\in A\cap T^{-m}A$ выполнено  неравенство
$$
  \langle  \J_x^{i(m,1)}\chi_B |\chi_B \rangle
     \ >\ \frac{1}{2q} \mu(B).
$$
\medskip

{ Доказательство} пункта (v).

Теперь   рассмотрим  оператор
$$\JJ_m:L_2(X,\mu)\to L_2(X,\mu)\otimes L_2(X,\mu),$$ определенный
формулой
$$
  \langle \JJ_m\chi_A |\chi_B\otimes\chi_C \rangle  =
  \int_A\langle \J_x^{i(m,1)}\chi_B |\chi_C \rangle d\mu(x).
$$
Для различных $m,k$ выполнено: или  пространства $\JJ_m L^2_0$, $\J_k L^2_0$
совпадают,
или  эти  пространства ортогональны (см. доказательство
теоремы 2.1.2).

Наша задача -- доказать, что множество попарно ортогональных пространств,
взятых из набора $\{\JJ_m L^2_0\}$, должно быть конечным.

Положим

$\bar{f}_m=\chi_{A\cap T^{-m}A} - \mu(A\cap T^{-m}A)\un$,

$\bar{\chi}_B=\chi_B -\mu(B)\un$.
\\
Так как выполнено
$$
     \langle \JJ_m\bar{f}_m |\un\otimes\bar{\chi}_B \rangle \ = \
     \langle \JJ_m\bar{f}_m |\bar{\chi}_B\otimes\un \rangle \ = \
    \langle \un |\bar{\chi}_B\otimes\bar{\chi}_B \rangle  \ =\  0,
$$
получаем
$$
    \langle \JJ_m\bar{f}_m |\bar{\chi}_B\otimes\bar{\chi}_B \rangle  =
  \langle \JJ_m\chi_{A\cap T^{-m}A} |\chi_B\otimes\chi_B \rangle -
  \mu(A\cap T^{-m}A)\mu(B)\mu(B).
$$
Таким образом, для всех больших $m$ выполнено неравенство
$$
    \langle \JJ_m\bar{f}_m |\bar{\chi}_B\otimes\bar{\chi}_B \rangle  \ > \
\frac{1}{2q}\mu(A)^2\mu(B) - \mu(A)^2\mu(B)^2 >0.
$$
Теперь предположим, что для бесконечного множества $N$ для  всех $m\in N$
 пространства $\{\JJ_m L^2_0\}$ попарно ортогональны.
Покажем, что предположение приводит к противоречию.

Действительно,
функции $\JJ_m\bar{f}_m$ попарно ортогональны, причем
для некоторого положительного числа $c_1$
выполнено   $$\|\JJ_m\bar{f}_m\|> c_1.$$
Но из этого вытекает следующее: для некоторой положительной константы $c_2$
для всех $m\in N$ выполнено неравенство
$$
    |\langle \JJ_m\bar{f}_m |\bar{\chi}_B\otimes\bar{\chi}_B \rangle|  \ > c_2,
$$
что невозможно, так как из-за попарной ортогональности $\JJ_m\bar{f}_m$
получается, что
$$\|\bar{\chi}_B\otimes\bar{\chi}_B\| =\infty.$$

Следовательно, найдется  множество $M$ положительной плотности
(не меньшей, чем число $\frac{1}{|N|}$)
такое, что  все пространства $\{\JJ_m L^2_0\}$ совпадают при $m\in M$.
Отсюда  следует, что  для некоторого $\J=\J_{m_0}$
$$
\forall m\in M \ \ \    \JJ_m =\J T^{k(m)}.
$$
Мы доказали, что  выполнено
$$
 \forall m\in M \ \ \ \ \ \
            \P_{T^m(x)}^\ast \P_x  = \frac{1}{p}(\J_{T^{k(m)}(x)}+\dots),
$$
причем  для всех $x\in A\cap T^{-m}A$
$$
  \langle  \J_{T^{k(m)}(x)}\chi_B |\chi_B \rangle
     \ >\ \frac{1}{2q} \mu(B).
$$
\\  

Чтобы завершить  доказательство утверждения 4.3, покажем необходимость
следующих условий:
 $$k(m)\to \infty, \ \ |k(m)-m|\to \infty.$$
Если для  бесконечного множества, элементы которого обозначим через $m'$,
выполнено  $k(m')=s$, где $s$ фиксировано, то получим
$$
       \Theta = \lim_{N\to\infty} \frac{1}{N}\sum_{i=1}^{N}
    \int \P_{T^{m'(i)}(x)}^\ast \P_x d\mu(x) =
       \frac{1}{p}(\J_{T^{s}(x)}+\dots).
$$
Следовательно, выполнено
$$\J_{T^{s}(x)}\equiv \Theta , \ \J_x \equiv \Theta.$$

Если   для  бесконечного множества различных $m'$
выполнено $k(m')-m'=s$, то
$$
       \Theta = \lim_{N\to\infty} \frac{1}{N}\sum_{i=1}^{N}
    \int \P_x^\ast \P_{T^{-m'(i)}(x)}d\mu(x)  = \frac{1}{p}(\J_{T^{s}(x)}+\dots).
$$
Снова получаем, что $\J_x \equiv \Theta$. Но это противоречит нетривиальности операторов $\J_x $.
Таким образом, $k(m)\to \infty$  и $|k(m)-m|\to \infty$.
 Этим завершается доказатеьство леммы 4.4, утверждения 4.3 и  теоремы 4.2. 

\vspace{4mm}
Для потоков    похожие задачи обнаружения   свойств $S_n$ решаются быстрее и эффектнее благодаря возможности использовать малые возмущения джойнингов  и  примененять  бесконечные башни расширений, см.  \cite{RT}.

\vspace{5mm}  Автор благодарит Жан-Поля Тувено за многочисленные 
интересные и полезные   обсуждения тематики статьи.

\end{fulltext}

\end{document}